\documentclass[journal ]{new-aiaa}
\usepackage[utf8]{inputenc}
\usepackage{textcomp}

\usepackage[dvipsnames]{xcolor}
\usepackage[normalem]{ulem}
\usepackage{mathtools}   
\usepackage{acronym}
\usepackage{tikz}
\usepackage{pgfplots}
\pgfplotsset{compat=1.18}
\usepackage{array}
\usepackage{hyperref}
\usepackage[capitalize]{cleveref}
\usepackage{subcaption}
\usepackage{placeins}
\usepackage{url}
\usepackage{float}
\usepackage{multirow}
\usepackage[symbols,nogroupskip,sort=standard,stylemods=mcols,style=mcolindex]{glossaries-extra}
\usepackage{chngcntr}
\usepackage{apptools}
\usetikzlibrary{positioning, arrows.meta, shapes.multipart} 
\usepackage{rotating} 

\usepackage{graphicx}
\usepackage{amsmath}
\usepackage[version=4]{mhchem}
\usepackage{siunitx}
\usepackage{longtable,tabularx}
\def\d{\mathrm{d}}

\def\p{\mathrm{p}}
\def\refr{\mathrm{ref}}
 
\def\s{\mathrm{s}}
\def\dT{\mathrm{T}}

\newcommand{\caos}{\mathrm{chaos}}

\newcommand{\plan}{\mathrm{plane}}

\newcommand{\pc}{\mathrm{PoC}}

\newcommand{\rad}{\mathrm{radial}}

\newcommand{\std}{\mathrm{std}}
\newcommand{\tca}{\mathrm{TCA}}

\newcommand{\cB}{\mathcal{B}}
\newcommand{\cC}{\mathcal{C}}

\newcommand{\cL}{\mathcal{L}}

\newcommand{\cS}{\mathcal{S}}
\newcommand{\cT}{\mathcal{T}}

\newtheorem{property}{Property}[section]

\newtheorem{definition}{Definition}[section]

\newtheorem{assumption}{Assumption}

\newcommand{\bi}{\mathbf{i}}
\newcommand{\bj}{\mathbf{j}}
\newcommand{\bk}{\mathbf{k}}

\newcommand{\bv}{\mathbf{v}}

\newcommand{\bx}{\mathbf{x}}

\newcommand{\Sig}{\Sigma}

\newcommand{\ox}{\overline{x}}

\newcommand{\te}{\mathrm{te.}}

\DeclareMathAlphabet{\mathmybb}{U}{bbold}{m}{n}

\newcounter{inlineenum}
\renewcommand{\theinlineenum}{\alph{inlineenum}}
\newenvironment{inlineenum}
{\unskip\ignorespaces\setcounter{inlineenum}{0}%
\renewcommand{\item}{\refstePoCounter{inlineenum}{\textit{\theinlineenum})~}}}
{\ignorespacesafterend}

    \acrodef{cnes}[CNES]{Centre National d'Études Spatiales}
    
    \acrodef{cdm}[CDM]{Collision Data Message}
    
    \acrodef{geo}[GEO]{Geosynchronous Orbit}
    
    \acrodef{hbr}[HBR]{Hard-Body Radius}
    \def\hbr{\ac{hbr}}
    \acrodef{isae}[ISAE-SUPAERO]{Institut Supérieur de l'Aéronautique et de l'Espace}
    
    \acrodef{leo}[LEO]{Low-Earth Orbit}
    
    \acrodef{poc}[PoC]{Probability of Collision}
    \def\poc{\ac{poc}}
    \acrodef{resp}[resp.]{respectively}
    \def\resp{\ac{resp}}
    \acrodef{spoc}[sPoC]{Scaled Probability of Collision}
    
    \acrodef{tca}[TCA]{Time of Closest Approch}
    \def\tca{\ac{tca}}

\begin{document}

\def\ourtitle{Probability of Collision with Tethered Spacecraft}

\def\NameYema{Yema \textsc{Paul}}

\def\DivSME{DOA/SME/SE}
\def\DivOR{DTN/DV/OR}
\def\OrgCNES{\acs{cnes}}
\def\StrCNES{18 Av. Édouard Belin}
\def\ZiPoCNES{31400}

\def\DivLab{Space Advanced Concepts Laboratory}
\def\OrgISAE{\acs{isae}}
\def\StrISAE{10 Av. Édouard Belin}
\def\ZipISAE{31400}

\def\CityToul{Toulouse}
\def\CountryFR{France}

\def\JobYema{PhD candidate, \DivLab}
\def\JobPau{Head of orbit determination office, \DivOR}

\author{
Yema Paul\(^{a,1}\) \quad
Emmanuel Delande\(^{b,2}\) \quad
Francois Vinet\(^{b,3}\) \quad
Francois Laporte\(^{b,4}\) \\[0.5ex]
Manuel Sanjurjo-Rivo\(^{c,5}\) \quad
Aldo Tonnini\(^{a,6}\) \quad
Joan-Pau Sanchez\(^{a,7}\)
}

\affil{\(^{a}\) ISAE-SUPAERO, 31400 Toulouse, France}
\affil{\(^{b}\) CNES, 31400 Toulouse, France}
\affil{\(^{c}\) Universidad Carlos III de Madrid, Leganés 28911, Spain}

\footnotetext[1]{PhD candidate, Department of Aerospace Vehicles Design and Control (DCAS), \texttt{yema.paul@isae-supaero.fr}}
\footnotetext[2]{Space surveillance specialist, \texttt{emmanuel.delande@cnes.fr}}
\footnotetext[3]{Space surveillance specialist, \texttt{francois.vinet@cnes.fr}}
\footnotetext[4]{Space surveillance specialist, \texttt{francois.laporte@cnes.fr}}
\footnotetext[5]{\texttt{msanjurj@ing.uc3m.es}}
\footnotetext[6]{\texttt{aldo.tonnini@student.isae-supaero.fr}}
\footnotetext[7]{Professor in Astrodynamics and Mission Design, \texttt{joan-pau.sanchez@isae-supaero.fr}}

\def\ourabstract{
}

\def\ourkeywords{
Conjunction analysis; Probability of collision; Tethered spacecraft; Collision risk estimation; Flexible tethers; Space situational awareness}

\title{\ourtitle}

\maketitle


\textbf{Keywords:~} \ourkeywords

\acresetall




\section{Introduction}

The increasing congestion of orbital space has led to a growing number of close encounters between operational spacecraft and other objects. In this context, the probability of collision ($\pc$) has become a central metric in conjunction assessment and decision-making \cite{alfano2022_CA_needs}. For most space objects, standard methods rely on modeling the spacecraft as a sphere with a predefined hard-body radius (HBR), and evaluating $\pc$ over that sphere under Gaussian uncertainty assumptions \cite{akella_probability_2000}.

However, certain classes of tethered spacecraft challenge this framework. In particular, systems composed of multiple bodies connected by a long, flexible tether present a spatially extended geometry and limited observability~\cite{cosmo1997tethers, chen2014tether_dynamics_review}. In many operational scenarios, only the main body is trackable, while the tether and its extremity remain unobserved. As a result, applying standard methods with  spherical HBR the length of the tether often yields overly conservative risk estimates—typically in the range of $10^{-2}$ to $1$—even for low-risk encounters. Moreover, the number of high-risk conjunctions involving tethered spacecraft is currently an operational preoccupation. It highlights the operational urgency of developing dedicated collision risk assessment methods that properly reflect the extended geometry and uncertain configuration of tethered systems.

Past studies have proposed analytical approaches that assume the tether configuration is known at the time of closest approach~\cite{patera2002method}. However, this assumption rarely holds in real scenarios. Tether dynamics are highly sensitive to environmental forces and initial conditions, and their shape is typically unobservable in catalog-based collision screening~\cite{chen2014tether_dynamics_review}.

In this work, we propose a new methodology to assess collision risk with tethered spacecraft in the presence of tether configuration uncertainty. Instead of relying on a known tether geometry, we compute the worst-case $\pc$ over the space of all physically admissible tether configurations. We consider three operationally relevant scenarios, each corresponding to a distinct constraint on the possible tether configurations:
\begin{itemize}
    \item \emph{Worst-case configuration}: the tether is assumed to explore any admissible 3D shape of fixed length, without constraint, representing a conservative upper bound.
    \item \emph{Planar}: the tether is constrained to lie within the orbital plane of the main vessel~\cite{misra2008tether_dynamics};
    \item \emph{Radial}: the tether remains aligned with the radial direction, as in Earth-pointing electrodynamic systems~\cite{li2024collision}.
\end{itemize}

This framework enables the derivation of conservative yet realistic upper bounds on collision risk with tethered systems. We apply the method to two real conjunction events involving a \SI{4}{\kilo\meter}-long inextensible flexible tethered spacecraft. Both events triggered operational alerts with $\pc$ values close to 1 using standard procedures. Our results offer a more accurate and less conservative risk estimate, supporting more informed decision-making for future tethered missions.

\section{Method}
\subsection{Definitions and generalities}
Consider the conjunction event between an inextensible tether spacecraft (so called primary object) and a secondary object.
We assume that some information is available on the states of the main body and the secondary object (e.g., drawn from a space catalog prior to the collision risk analysis).
The tether spacecraft consists in a tether of length $\ell$ linking the main vessel -- the \emph{main body} -- to a mass -- the \emph{small body}. By extension, the main body, small body, and secondary objects are called \emph{(rigid) bodies}. The \emph{(kinematic) state} of the bodies describe their position and velocity coordinates in some reference inertial frame. We set the encounter problem in the usual frame of work of fast encounters. The period around the \tca{} during which the collision may occur is the \emph{encounter window}, and we assume the following:

\begin{assumption}[Dynamics of the main body and secondary object] \label{hyp:object_dynamics}
Throughout the encounter window:
    \begin{enumerate}
        \item The main body (\resp{}, the secondary object) evolves on a straight line, with constant velocity.
        \item The velocity of the main body (\resp{} the secondary object) is known, and given by $\mathbf{v}_{\p_1}$ (\resp{} $\mathbf{v}_{\s}$).
    \end{enumerate}
\end{assumption}




\begin{definition}[Encounter frame and conjunction plane]
Let $\bx_{\p_1}$ be the nominal position of the main body at time of closest approach.  
We consider the orthonormal frame $(\bi, \bj, \bk)$ given by \cite{akella2000probability}:
    \begin{align}\label{eq:ijk conjunction plane}
            \bi = \frac{\bv_{\s}^{}-\bv_{\p_1}^{}}{|\bv_{\s}^{}-\bv_{\p_1}^{}|}, \quad
        \mathbf{j} = \frac{\bv_{\s} \wedge \bv_{\p_1}}{|\bv_{\s} \wedge \bv_{\p_1}|}, \quad \bk = \bi \wedge \bj.
    \end{align}
We define the \emph{local encounter frame} as $\cL = (\bx_{\p_1}; \bi, \bj, \bk)$, and 
the \emph{conjunction plane} $\cC$ as the 2D affine plane through $\bx_{\p_1}$ spanned by $(\bj, \bk)$:
\begin{align}
    \cC = \bx_{\p_1} + \mathrm{span}(\bj, \bk).
\end{align}
\end{definition}

\begin{assumption}[Information on the main body and secondary object]\label{hyp:object_info}
Throughout the encounter window, the position uncertainty of the main body (\resp{}, the secondary object) is described by a Gaussian distribution with constant covariance.
\end{assumption}

At \tca{}, the projection of the nominal (or mean) value of the position of the main body (\resp{}, the secondary object) in the conjunction plane is denoted by $x_{\p_1}$ (\resp{} $x_{\s}$). Likewise, the projection of the (constant) covariance on the position of the main body (\resp{}, the secondary object) in the conjunction plane is denoted by $\Sigma_{\p_1}$ (\resp{} $\Sigma_{\s}$).

Given the geometries of the bodies involved in the collision event, the tether is assimilated to a line of negligible thickness. We impose little physical constraints to its shape and mechanical properties, but for its fixed length $\ell$. In order to describe the shape of the tether throughout the encounter window, we make use of the \emph{local encounter frame} $\cL$ with the tuple $(\bi, \bj, \bk)$ centered on the main body's nominal position. Then:

\begin{definition}[Tether configuration]
    \label{def:tether_shape}
    A \emph{(tether) configuration} is a continuous mapping $\bx_{\te}:~[0, 1] \mapsto \mathbb{R}^3$ such that
    \begin{align}\label{eq:tether_position}
        \bx_{\te}(0) = 0, \quad \mathrm{and} \int_0^1 \left\| \frac{\d \bx_{\te}(s)}{\d s} \right\| \d s &= \ell,
    \end{align}
    describing the shape of the tether in the local encounter frame  \(\mathcal{L}\). In particular, the tip $\bx_{\te}(0)$ coincides with the nominal position of the main body, and the tip $\bx_{\te}(1)$ coincides with the center of gravity of the small body. The set of all configurations is denoted by $\cS_{\te}$. 
\end{definition}

We also denote by $R_{\p_1}$ (\resp{} $R_{\p_2}$, $R_{\s}$) the \hbr{} of the main body (\resp{} small body, secondary object). Finally, we denote by $(\mathbf{R}_{\p}, \mathbf{T}_{\p}, \mathbf{N}_{\p})$ the local RTN frame of the main body, where $\mathbf{R}_{\p}$ (\resp{} $\mathbf{T}_{\p}$, $\mathbf{N}_{\p}$) is the radial (\resp{} tangential, out-of-plane) component.

\subsection{Problem Formulation}
A detailed modeling of tether dynamics is beyond the scope of this paper. This is further motivated by the fact that (1) the tethered system can exhibit highly complex behavior in orbit, including bounded, unbounded, or even chaotic configurations under certain energy and orbit conditions~\cite{misra2008tether_dynamics}; and (2) neither the tether nor the small body is typically observable with sufficient accuracy using ground-based sensors. For the sake of this problem, we will assume that

\begin{assumption}[Tether dynamics]
    \label{hyp:tether_dynamics}
    Throughout the collision encounter, the configuration of the tether is fixed and characterized by some unknown distribution $\rho_{\te}$ on $\cS_{\te}$.
\end{assumption}

\begin{figure}
    \centering
    \includegraphics[width=0.7\linewidth]{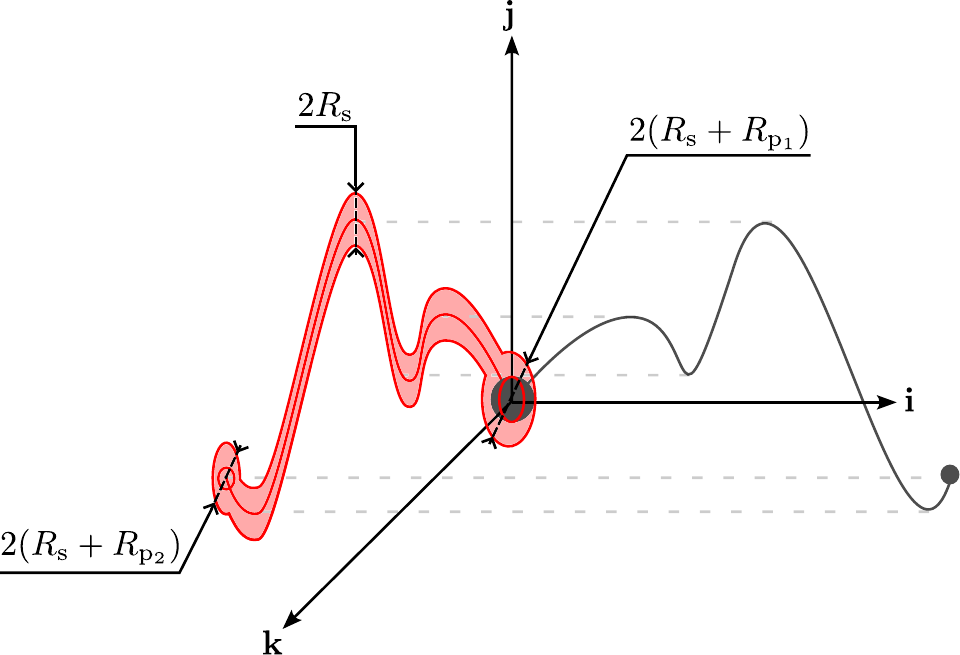}
    \caption{\small Projection of the primary object in the conjunction plane, with the collision region depicted in red\label{fig_proj_cplane}}
\end{figure}

Further denote by $\varphi_{\cL}^{\cC}$ the projection from the local encounter frame $\cL$ to the conjunction plane $\cC$. Then, we can formulate the collision risk assessment problem as a generalization of the classical, 2D conjunction problem \cite{akella2000probability}. Indeed, the geometry of the primary object no longer reduces to a sphere with radius $R_{\p_1}$ centered on the main body's nominal position $x_{\p_1}$: it also incorporates a sphere with radius $R_{\p_2}$ accounting for the small body, as well as the projection $x_{\te} \coloneq \varphi_{\cL}^{\cC} \circ \bx_{\te}$ representing the "footprint" of the tether in the conjunction plane (see \cref{fig_proj_cplane}).
We define the \textit{Combined Tether Hard Body} by 
\begin{align}
    \cT(x_{\te}; R_{\p_1}, R_{\p_2}, R_{\s}) \coloneq 
    \cB_2(x_{\te}(0), R_{\p_1}+R_{\s}) \cup
    \bigcup_{0\leq t \leq 1} \cB_2(x_{\te}(t), R_{\s})
    \cup \cB_2(x_{\te}(1), R_{\p_2}+R_{\s}).
\end{align}
where $\cB_2$ denotes 2-dimensional balls in the conjunction plane $\cC$, and $\bigcup_{0\leq t\leq 1} \cB_2(x_{\te}(t), R_{\s}) $ is the union of the 2-dimensional balls sweeping alongside the projection of the tether in the conjunction plane.

The PoC \cite{akella2000probability, hall2019_PoC_2D_limitations} between the tether spacecraft and the secondary is then
\begin{align} \label{eq_pcref}
    \pc_{\refr} =
    &  \int  
    \pc\Big( \cT(x_{\te}; R_{\p_1}, R_{\p_2}, R_{\s}); ~ x_{\s}, \Sig_{\p_1}+\Sig_{\s} \Big)\rho_{\te}(\bx_{\te}) \d\bx_{\te} ,
\end{align}
where\begin{align}
    \pc(A; \ox, \Sigma)=\frac{1}{2\pi \sqrt{|\det \Sigma|}}  \int_A 
    \exp\left(
        -\frac{1}{2}(x-\ox)^T \Sigma^{-1}(x-\ox)
    \right)\d x.
\end{align}
Unfortunately, $\rho_{\te}$ is typically unknown to the conjunction analyst, so the computation of the reference value Eq.\eqref{eq_pcref} remains out of reach in most practical situations.

\subsection{Evaluation of collision risk}

We will then consider particular cases, corresponding to various physical configurations of the tether, and draw the relevant information on the \poc{}.

\subsubsection{Maximum PoC over All Tether Configurations}

We first consider a conservative upper bound on the reference value $\pc_{\refr}$ in Eq.~\eqref{eq_pcref}, by maximizing the collision probability over all physically admissible tether configurations. Since neither the tether nor the small body is tracked, and no further dynamical assumptions are made, all configurations in $\cS_{\te}$ are considered possible.
Although tethered systems can exhibit bounded, unbounded, or even chaotic behavior~\cite{misra2008tether_dynamics}, this formulation simply reflects a lack of knowledge. We refer to the resulting worst-case estimate as the $\overline{\pc}_{\caos}$ defined by
\begin{align}\label{eq:maximum_over_configuration}
    \overline{\pc}_{\caos} = \sup_{\bx_{\te}\in \cS_{\te}}
        \pc\Big( \cT(x_{\te}; R_{\p_1}, R_{\p_2}, R_{\s}); ~ x_{\s}, \Sig_{\p_1}+\Sig_{\s} \Big).
\end{align}
The formal proof that $\overline{\pc}_{\caos}$ is an upper bound on $\pc_{\refr}$ is provided in \cref{sec:proof_PoC_caos_upper_bound}.

\subsubsection{Maximum PoC over all Planar Configurations}

An important case is when 
the small body, the tether and the main body of the tether spacecraft remain co-planar \cite{misra2008tether_dynamics},
either naturally by the dynamics or by active control.
\begin{assumption}[Tether Planar Dynamics]
    \label{hyp:planar_dynamics}
    The small body and the tether remain in the same orbital plane as the main body. Otherwise said, $\bx_{\te}(t)\cdot \mathbf{N}_{\p} = 0$, for all $0 \leq t \leq 1$. 
\end{assumption}
Define the matrices
\begin{equation}
    P =
    \begin{bmatrix}
        \mathbf{R}_{\p} \cdot \mathbf{j} & \mathbf{T}_{\p} \cdot \mathbf{j} \\
        \mathbf{R}_{\p} \cdot \mathbf{k} & \mathbf{T}_{\p} \cdot \mathbf{k}
    \end{bmatrix}
    \text{ and } Q = \ell^2 PP^{\dT}.
\end{equation}
\cref{hyp:planar_dynamics} on tether dynamics  constrains the tethered spacecraft system $\bx_{\te}$ to remain 
in a 2-dimensional ball in the plane $(\mathbf{R}_{\p}, \mathbf{T}_{\p})$, whose projection onto the conjunction plane $\cC$ is a 2-dimensional ellipse $E_2(x_{\p_1}, Q)$ defined by
the implicit equation   
\begin{align}
    x_{\te}^{\dT}~Q^{-1} x_{\te}\leq 1.
\end{align}
Therefore, by considering the set of planar configurations 
\begin{align}
    \cS_{\plan}\coloneq \{\bx_{\te} \in \cS_{\te}~|~x_{\te}\subset E_2(x_{\p_1}, Q)\},
\end{align}
\cref{eq:maximum_over_configuration} becomes to 
\begin{align}
    \overline{\pc}_{\plan} = \sup_{\bx_{\te}\in \cS_{\plan}}
        \pc\Big( \cT(x_{\te}; R_{\p_1}, R_{\p_2}, R_{\s}); ~ x_{\s}, \Sig_{\p_1}+\Sig_{\s} \Big) 
\end{align}

\subsubsection{PoC for straight configuration}
Finally, we consider the specific case where the tether is fully stretched and aligned with Earth's gravitational pull on the main body, a configuration that may occur in particular with gravity-gradient stabilized or electrodynamic tethers \cite{cosmo1997tethers}.
\begin{assumption}[Straight Earth-Pointing Tether]
    \label{hyp:straight_earth_pointing}
    The tether is fully deployed and modeled as a straight line aligned with the local radial direction pointing toward Earth. That is, the configuration of the tether in the local encounter frame $\cL$ is
    \begin{align}
    \bx_{\te, \rad}(t) = -t \, \ell \, \mathbf{R}_{\p}, \quad \text{for all } t \in [0,1].
    \end{align}
    In particular, the tether lies entirely within the orbital plane and experiences no transverse motion or deformation.
\end{assumption}
Since the tether lies along $\mathbf{R}_{\p}$, its projection onto the conjunction plane $\cC$ is given by
\begin{align}
    x_{\te, \rad}(t) = - t\ell \, u, \quad \text{where} \quad
u = 
\begin{bmatrix}
\mathbf{R}_{\p} \cdot \mathbf{j} \\
\mathbf{R}_{\p} \cdot \mathbf{k}
\end{bmatrix},
\end{align}
for all  $t \in [0,1]$. Thus, \cref{eq:maximum_over_configuration} becomes to 
\begin{align}
    \pc_{\rad} =
        \pc\Big( \cT(x_{\te, \rad}; R_{\p_1}, R_{\p_2}, R_{\s}); ~ x_{\s},
                    \Sig_{\p_1}+\Sig_{\s} \Big), 
\end{align}
that is, the analyst has now access to the value of the PoC.

\subsection{Comparison of PoC methods}
We analyze in this section the several methods presented above. We first introduce the standard operational method that spares the modeling of the tether altogether and accounts for its existence solely through an inflated \hbr{} for the primary object, namely, consider \hbr{} $\ell$ instead of $R_{\p_1}$. The corresponding PoC is then
\begin{align} \label{eq_poc_std}
    \pc_{\std} = \pc \Big(
        \cB_2(x_{\p_1}, \ell+R_{\s})
    ; ~ x_{\s}, \Sig_{\p_1}+\Sig_{\s}\Big).
\end{align}

This "standard" method is mostly advantageous due to its simplicity, since it avoids the modeling of the tether and remains within the remit of the classical 2D PoC approach \cite{akella2000probability}. It does, however, correspond to a non-physical configuration of the tether that would \emph{simultaneously fill up} the entire 3-dimensional ball of radius $\ell$ centered on the main body, at any given time throughout the encounter. Thus it (largely) overestimates the collision risk, in the general case. More precisely, it holds that:
\begin{property}[Comparison of PoC methods]
    Under chaotic motion (i.e., no further assumptions beyond \cref{hyp:tether_dynamics}),
    \begin{align}
        \pc_{\refr} \leq \overline{\pc}_{\caos} \leq 
        \pc_{\std}.
    \end{align}
    Under planar motion (see \cref{hyp:planar_dynamics}),
    \begin{align}
        \pc_{\refr} \leq \overline{\pc}_{\plan} \leq 
        \overline{\pc}_{\caos} \leq 
        \pc_{\std}.
    \end{align}
    Under Earth-pointing configuration (see \cref{hyp:straight_earth_pointing}),
     \begin{align}
        \pc_{\refr} = \pc_{\rad} \leq \overline{\pc}_{\plan} \leq 
        \overline{\pc}_{\caos} \leq 
        \pc_{\std}.
    \end{align}
\end{property}


\section{Results and Discussion}




We illustrate the computation of the \poc{} using 
two real conjunction events involving a \SI{4}{\kilo\meter}-long, inextensible 
flexible tethered spacecraft (see \cref{tab:tether_caracteristics}).  
These events are referred to as \emph{Case 1} and \emph{Case 2}. The standard method in \cref{eq_poc_std} yields operational collision probabilities $\pc_{\std}$ 
of \SI{98}{\percent} and \SI{0.33}{\percent}, 
respectively, at \SI{15}{\hour} and \SI{38}{\hour} before the \tca{}.
\begin{table}[htbp]
    \centering
    \small
    \renewcommand{\arraystretch}{1.2}
    \begin{tabular}{ccccc}
        \hline
        Tether Type  & $\ell$ [$\si{\meter}$] & $R_{\p_1}$ [$\si{\meter}$]  & $R_{\p_2}$ [$\si{\meter}$] & $R_{\s}$ [$\si{\meter}$]  \\
        \hline
        Inextensible flexible& 4000 &     5  &   1    &  1 \\
        \hline
    \end{tabular}
    \caption{\small Tethered Spacecraft and Secondary Characteristics}
    \label{tab:tether_caracteristics}
\end{table}

\begin{figure}[htbp]
    \centering
    \begin{subfigure}[t]{0.45\textwidth}
        \includegraphics[width=\linewidth]{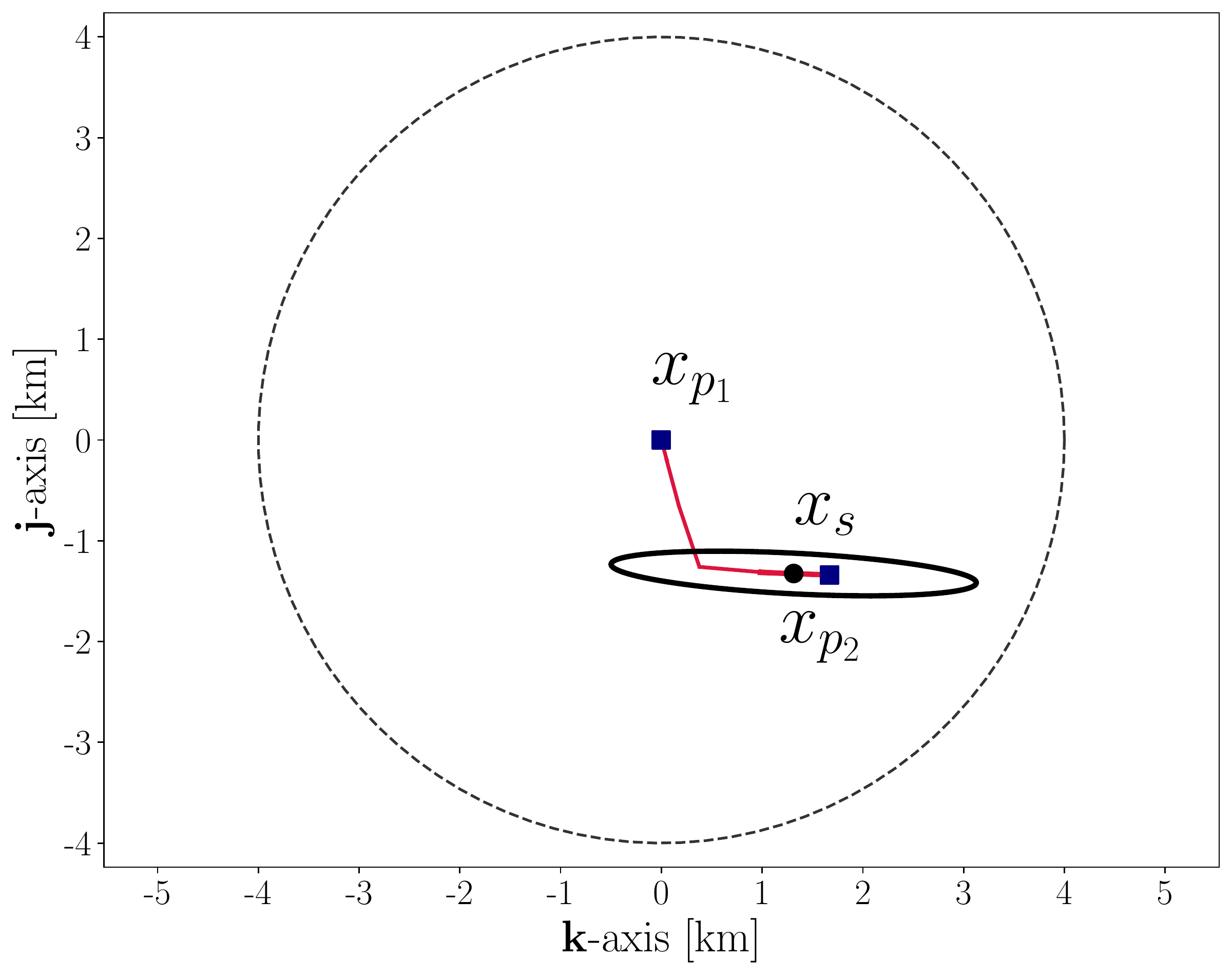}
        \caption{\small Case 1 under Chaotic Tether Motion}
        \label{fig:case1_chaos}
    \end{subfigure}
    \hfill
    \begin{subfigure}[t]{0.45\textwidth}
        \includegraphics[width=\linewidth]{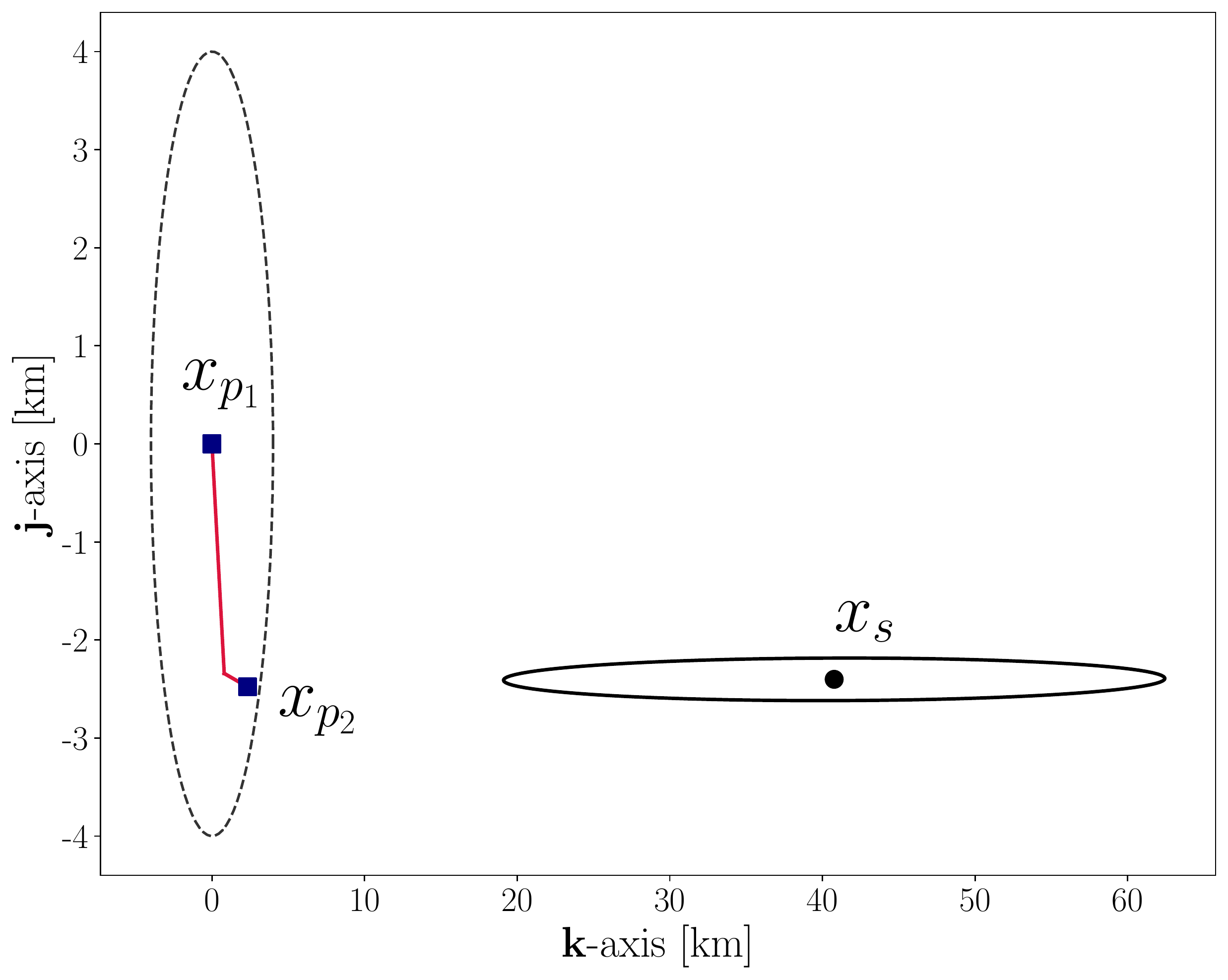}
        \caption{\small Case 2 under Chaotic Tether Motion}
        \label{fig:case2_chaos}
    \end{subfigure}

    \vspace{0.8em}

    \begin{subfigure}[t]{0.45\textwidth}
        \includegraphics[width=\linewidth]{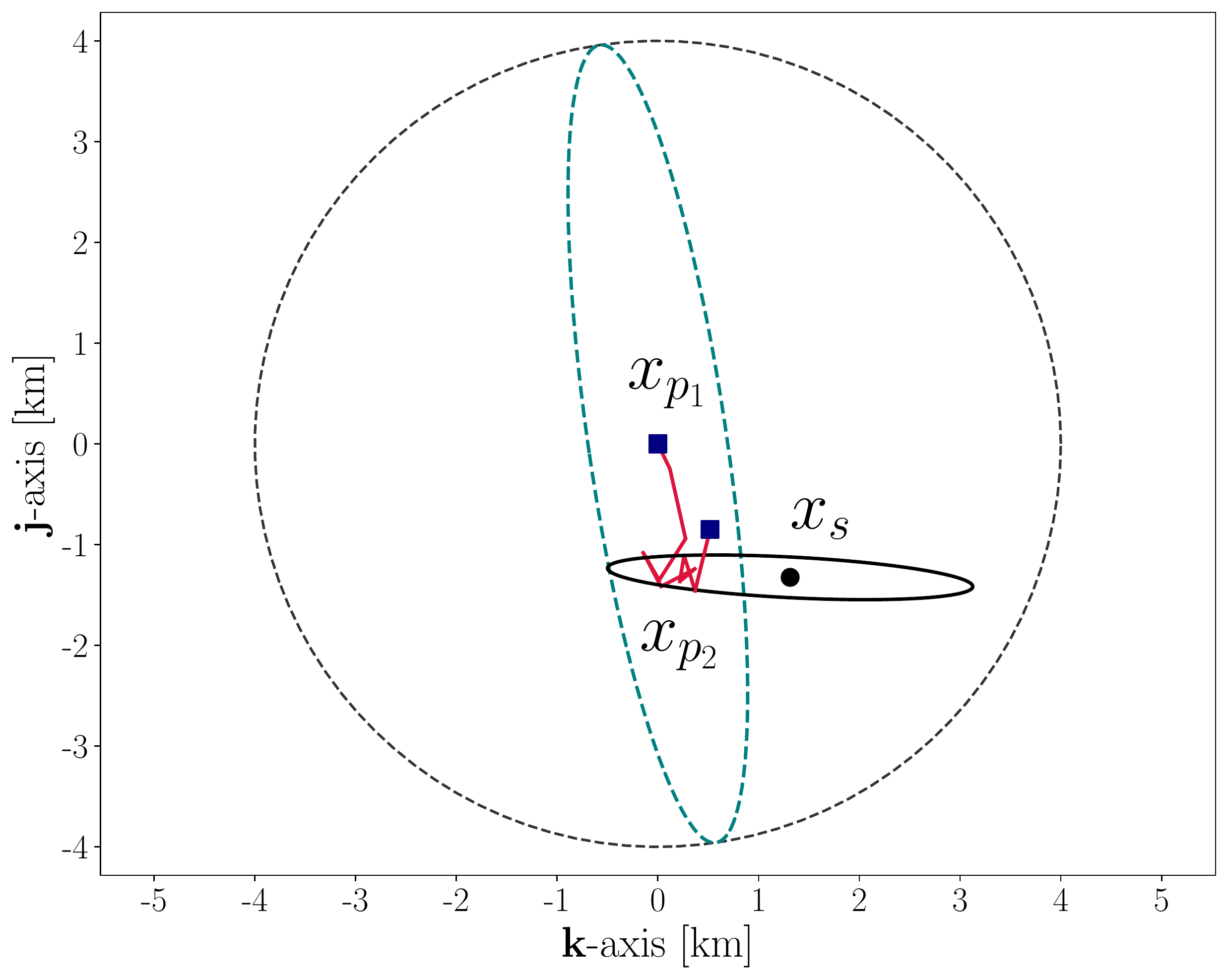}
        \caption{\small Case 1 under Planar Tether Motion}
        \label{fig:case1_planar}
    \end{subfigure}
    \hfill
    \begin{subfigure}[t]{0.45\textwidth}
        \includegraphics[width=\linewidth]{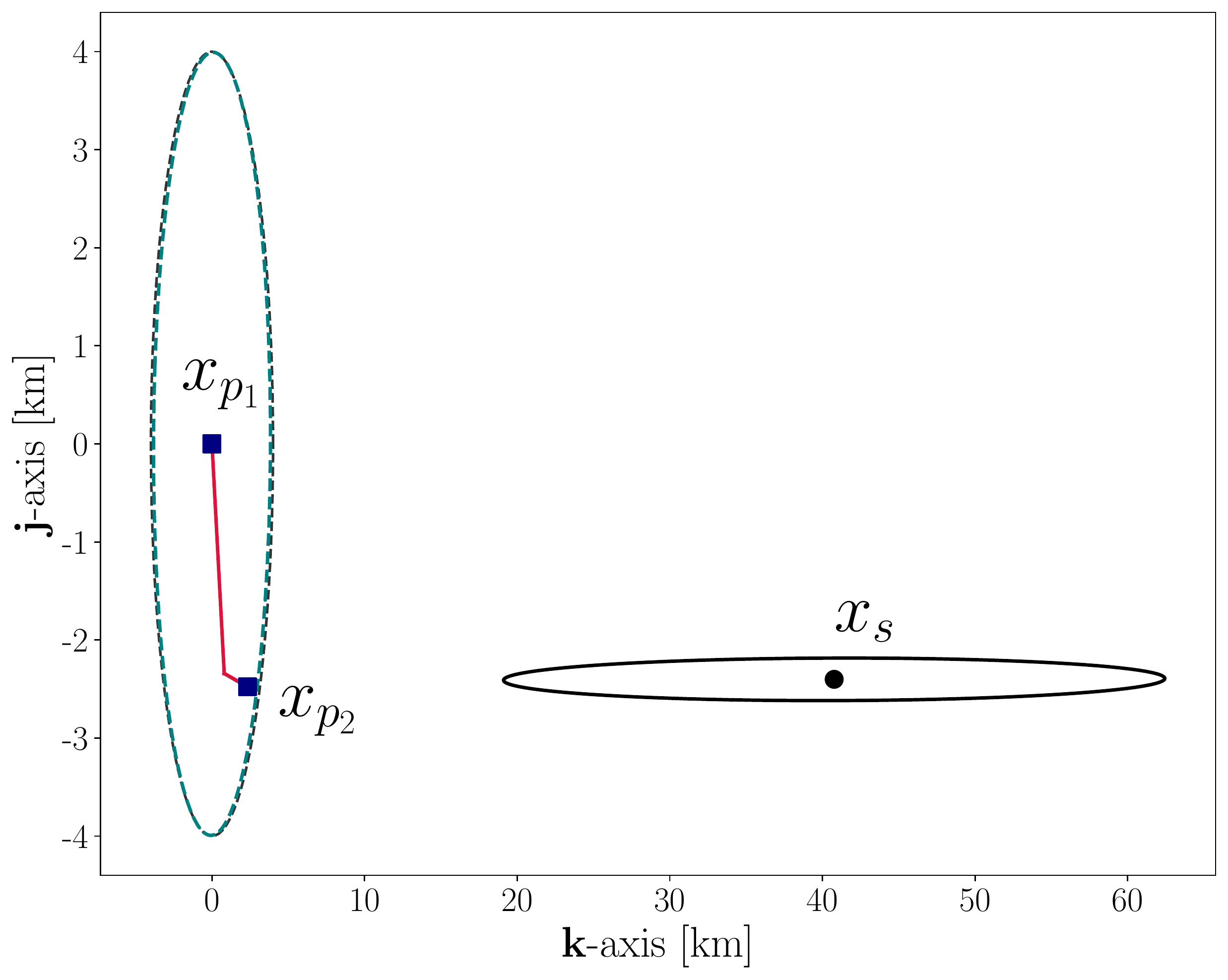}
        \caption{Case 2 under Planar Tether Motion}
        \label{fig:case2_planar}
    \end{subfigure}

    \vspace{0.8em}

    \begin{subfigure}[t]{0.45\textwidth}
        \includegraphics[width=\linewidth]{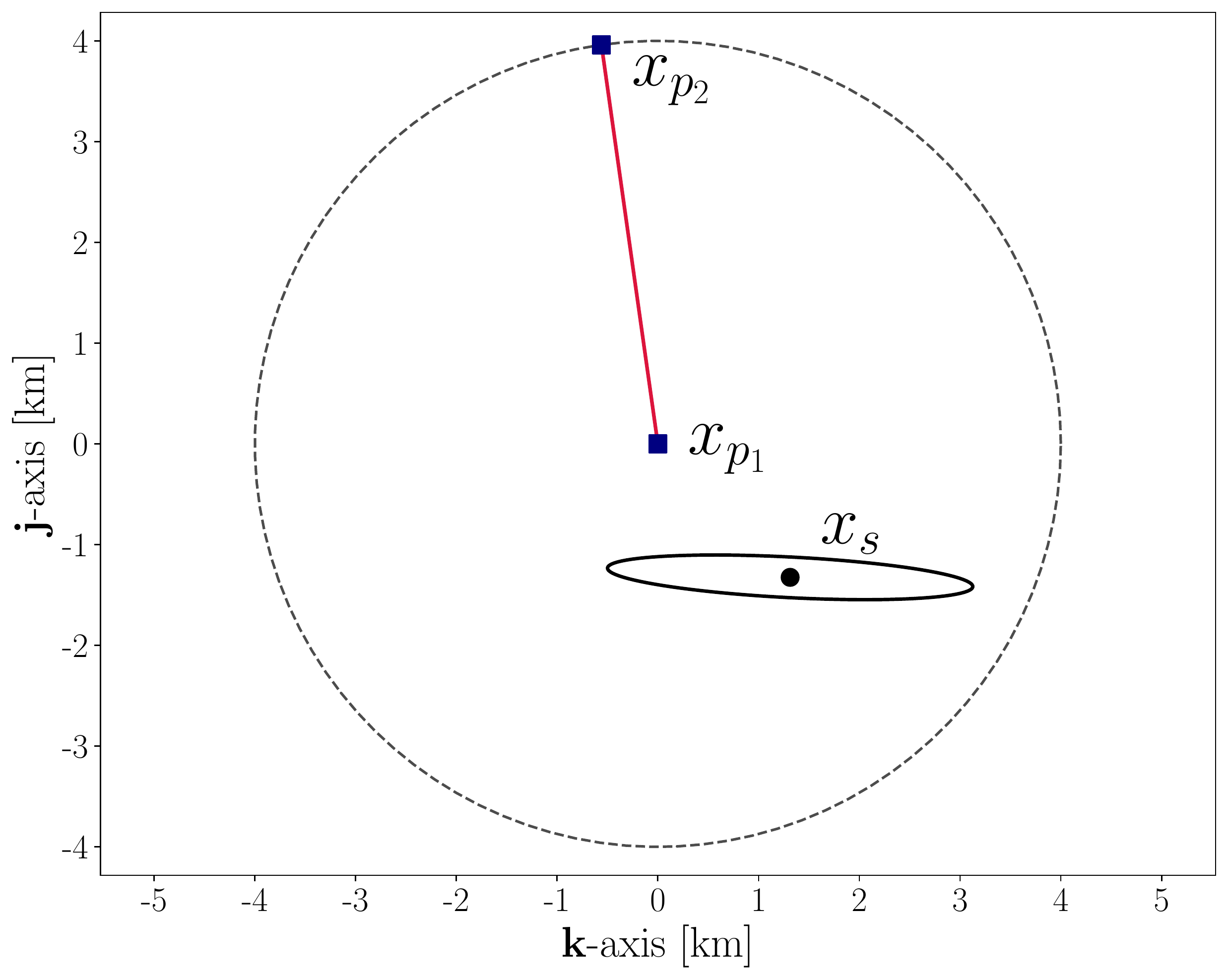}
        \caption{\small Case 1 with Earth Pointing Tether}
        \label{fig:case1_radial}
    \end{subfigure}
    \hfill
    \begin{subfigure}[t]{0.45\textwidth}
        \includegraphics[width=\linewidth]{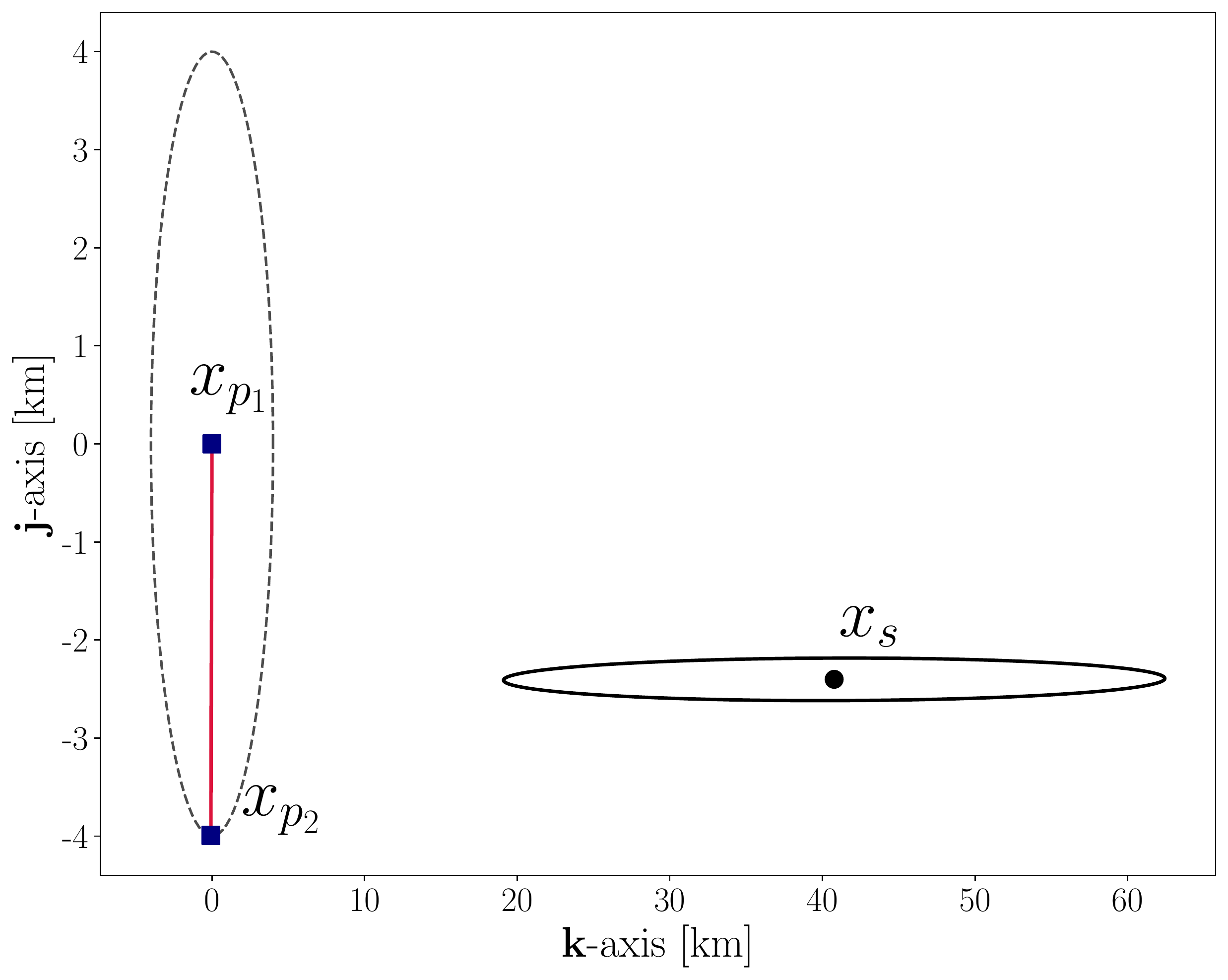}
        \caption{\small Case 2 with Earth Pointing Tether}
        \label{fig:case2_radial}
    \end{subfigure}

    \caption{\small Tether Configurations yielding the highest \poc{}, under Chaotic, Co-Planar, and Radial Motion Assumptions}
    \label{fig:tether_grid}
\end{figure}

To enable consistent comparison between scenarios, we assume a fixed \hbr{} of \SI{1}{\meter} for the secondary object, as specified in \cref{tab:tether_caracteristics}. In \cref{fig:tether_grid}, the tether is shown in red, with the main and small bodies represented by navy blue squares at $x_{\p_1}$ and $x_{\p_2}$, respectively. The $1\sigma$ confidence ellipse of the combined positional covariance, $\Sig_{\p_1} + \Sig_{\s}$, is centered on the secondary position $x_{\s}$. The black dashed circle represents the operational \hbr{} disk used in the standard approach.


\begin{table}[htbp]
    \centering
    \small
\renewcommand{\arraystretch}{1.6} 
    \begin{tabular}{lcccccc}
        \hline
         & Time to \tca{} & $\pc_{\std}$ & miss distance [$\si{\meter}$] 
         & relative speed [$\si{\kilo\meter\per\second}$] 
         & $ x_\s-x_{\p_1} ~~[\si{\meter}]$ & $\Sigma_{\p_1}+\Sigma_{\s} ~~[$\si{\meter\squared}$]$ 
         \\
        \hline
        Case 1 &
        15h & $98\%$ & 1865 & 14.808 & 
        $\begin{bmatrix} 
        1325 \\ 1313 
        \end{bmatrix}$ & 
        $\begin{bmatrix} 
            21446 & -72258 \\ 
            -72258 & 1435168 
        \end{bmatrix}$
        \\
        \hline
        Case 2 &
        38h & $3.3\times 10^{-3}$ & 40841 & 4.603 & 
        $\begin{bmatrix} 
        -2401 \\ 40770 
        \end{bmatrix}$
        & 
        $\begin{bmatrix} 
            20569 & 85842 \\ 
            85842 & 205011218 
        \end{bmatrix}$ \\
        \hline
    \end{tabular}
    \caption{\small Conjunction Variables at \tca{}.}
    \label{tab:conjunction_variables}
\end{table}

The primary conjunction variables at \tca{} are summarized in \cref{tab:conjunction_variables}. The conjunction-plane projections of tether configurations are illustrated for Case~1 in \cref{fig:case1_chaos,fig:case1_planar,fig:case1_radial} and for Case~2 in \cref{fig:case2_chaos,fig:case2_planar,fig:case2_radial}. The key difference between these two scenarios lies in the relative miss-distance: in Case~1, it is smaller than the tether's length, while in Case~2, it is significantly larger (see again \cref{tab:conjunction_variables}).

In the absence of dynamical constraints on the tether, all configurations within a 2D ball of radius $\ell$ are admissible. In order to estimate $\overline{\pc}_{\caos}$ and $\overline{\pc}_{\plan}$, an optimization algorithm was developed to explore tether configurations maximizing \poc{} values, assuming they can be approximated as a chain of \( N \) rigid bars connected in sequence (see Appendix \ref{app_algo}). 
This flexible representation also accommodates the special case of a rigid tethered spacecraft by setting $N=1$, thereby reducing the model to a single fixed segment between the two bodies.
The configurations shown in \cref{fig:case1_chaos,fig:case2_chaos} approximate the worst-case configurations yielding the upper bound $\overline{\pc}_{\caos}$. Under the planar motion  in \cref{hyp:planar_dynamics}, the approximates of the worst-case configurations are depicted in \cref{fig:case1_planar,fig:case2_planar}, and the corresponding admissible regions $\cS_{\plan}$ are bounded by green dashed curves. For Earth-pointing tethered spacecraft, the unique valid configuration is illustrated in \cref{fig:case1_radial,fig:case2_radial}. Finally, the collision risk estimates are gathered in \cref{tab:PoC_values}.


\begin{table}[htbp]
    \centering
    \small
\renewcommand{\arraystretch}{1.6} 
    \begin{tabular}{lcccccc}
        \hline
         & Time to TCA & $\pc_{\std}$ & $\overline{\pc}_{\caos}$ & $\overline{\pc}_{\plan}$ & $\pc_{\rad}$ \\
        \hline
        Case 1 &
        15h            & $98\%$       & $5.3\times 10^{-3}$      & $2.8\times 10^{-3}$ &      $3.7\times 10^{-24}$       \\
        \hline
        Case 2 &
        38h &    $3.3\times 10^{-3}$ & $5.8\times 10^{-6}$       & $5.8\times 10^{-6}$ & $9.6\times 10^{-7}$ \\
        \hline
    \end{tabular}
    \caption{\small Tether PoC for Two Conjunction Events.}
    \label{tab:PoC_values}
\end{table}
In both scenarios, the standard \poc{} estimate (i.e., $\pc_{\std}$) far exceeds typical operational thresholds, which generally lie between $10^{-5}$ and $10^{-4}$ \cite{alfano2022_CA_needs, florian2024future}. However, in the case of extended objects like tethered spacecraft, this method is considered to be overly conservative. As a result, even a probability estimate close to $1$ does not necessarily trigger a collision avoidance maneuver. This disconnect highlights a critical issue: when $\pc_{\std}$ fails to discriminate between high-risk and low-risk situations, its operational relevance is undermined—especially in scenarios where the miss distance between the secondary and the main body falls within the tether's physical range.

The refined estimates $\overline{\pc}_{\caos}$, $\overline{\pc}_{\plan}$, and $\pc_{\rad}$ for both scenarios are reported in \cref{tab:PoC_values}. 
Without any assumptions on the tether dynamics, $\overline{\pc}_{\caos}$ can be taken as the most conservative estimate.
 These values provide a more informative and nuanced characterization of the conjunction risk.   While the value of $\overline{\pc}_{\caos}$ remains high for Case~1, indicating a potentially concerning configuration, it is significantly lower for Case~2, suggesting a low-risk situation that would not have led to the overly pessimistic assessment of the standard method.


\subsection*{Conclusion}

The proposed methodology enables a more nuanced assessment of collision risk for tethered spacecraft by accounting for their spatial reach without relying on overly conservative assumptions. The three estimates $\overline{\pc}_{\caos}$, $\overline{\pc}_{\plan}$, and $\pc_{\rad}$ provide a 
spectrum of risk assessments, depending on the degree of knowledge or assumptions regarding the tether's dynamics. 
Notably, the comparison with the standard $\pc_{\std}$ estimate highlights the inadequacy of standard \poc{} calculation for tethered spacecraft. 

This approach offers a practical and interpretable framework for characterizing the spatial extent of tethered systems in conjunction analysis, and can be readily adapted to other elongated objects with a dominant dimension.
Instead of approximating the system as a sphere with radius equal to the tether 
length, it considers the set of physically admissible configurations—leading to a less conservative and more realistic assessment of collision risk. 
While final decisions on collision avoidance maneuvers rest with operational teams, 
the proposed methodology enhances existing practices by distinguishing between cases 
that are truly capable of resulting in a high probability of collision and those that,
although flagged as high-risk by standard methods, do not admit any physically 
plausible high-risk scenario. This distinction helps reduce unnecessary maneuvers 
without compromising safety.



\bibliography{biblio/references}   

\appendix

\section{Proof that PoC under Chaos Assumption is an Upper Bound on the Reference PoC}
\label{sec:proof_PoC_caos_upper_bound}
We show that the $\pc_{\caos}$, provides a conservative upper bound on the reference probability of collision $\pc_{\refr}$ defined in Eq.~\eqref{eq_pcref}. By definition:
\begin{align}
    \pc_{\refr} &=
    \int  
    \pc\Big( \cT(\bx_{\te}; R_{\p_1}, R_{\p_2}, R_{\s}); ~ x_{\s}, \Sigma_{\p_1} + \Sigma_{\s} \Big)
    \rho_{\te}(\bx_{\te}) ~\d\bx_{\te} \\
    &\leq
    \int  
    \sup_{\bx_{\te} \in \cS_{\te}}
    \pc\Big( \cT(\bx_{\te}; R_{\p_1}, R_{\p_2}, R_{\s}); ~ x_{\s}, \Sigma_{\p_1} + \Sigma_{\s} \Big)
    \rho_{\te}(\bx_{\te}) ~\d\bx_{\te} \\
    &=
    \pc_{\caos}
    \int  
    \rho_{\te}(\bx_{\te}) ~\d\bx_{\te} \\
    &= 
    \pc_{\caos}
\end{align}
where the last equality holds since $\rho_{\te}$ is a probability density over $\cS_{\te}$ and integrates to 1.

\section{Exploration of Tether Configurations Maximizing PoC Values} \label{app_algo}

This section provides an outline of the algorithm optimizing the tether configuration to maximize \poc{} values.
\subsection*{Modeling of Inextensible Tether Using \texorpdfstring{$N$}{N}-Rigid Bars}

We model the inextensible flexible tether as a chain of \( N \) rigid bars connected in sequence. Each bar is defined by its length \( \ell_k > 0 \) and orientation angle \( \theta_k \in [0, 2\pi) \) in the 2D conjunction plane \( \mathcal{C} \), with the constraint
\begin{align}
    \sum_{k=1}^{N} \ell_k \leq \ell,
\end{align}
where \( \ell \) is the total tether length.

Starting from the main body at the origin \( x_{\p_1} = 0 \), the configuration of the tether is defined recursively as
\begin{align}
    x_{\te}(k/N) = x_{\te}((k-1)/N) + \ell_k
    \begin{bmatrix}
        \cos \theta_k \\
        \sin \theta_k
    \end{bmatrix}, \quad \text{for } k = 1, \dots, N,
\end{align}
where the endpoint \( x_{\te}(1) =x_{\p_2}\) corresponds to the small-body location. 

Under planar motion model (see Assumption \ref{hyp:planar_dynamics}) the tether remains within a valid configuration space (e.g., resulting from projection of the 3D orbital plane into the conjunction plane), and we impose an elliptical feasibility constraint:
\begin{align}
    x_{\te}(k)^\top Q^{-1} x_{\te}(k) \leq 1, \quad \forall k = 1, \dots, N,
\end{align}
where \( Q \in \mathbb{R}^{2 \times 2} \) defines the shape and orientation of the allowable region in the conjunction plane \( \mathcal{C} \).

\subsection*{Optimization of Worst-Case Configuration}

To assess the worst-case configuration of the tether in a conjunction scenario, we aim to maximize the \poc{} induced by the tether geometry under the aforementioned constraints. The optimization variables are the number of segments, the segment lengths and angles:
\begin{align}
N \geq 1, \quad
\{ \ell_k, \theta_k \}_{k=1}^N, \quad \ell_k > 0, \ \theta_k \in [0, 2\pi),
\end{align}
subject to
\begin{align}
    \sum_{k=1}^{N} \ell_k \leq \ell.
\end{align}
The objective is to maximize the total collision probability contribution from 
main-body, tether constituted by $N$-rigid bars, and small-body.
However, in order to fasten the computation time, instead of computing 
\begin{align}
        \pc\Big( \cT(x_{\te}; R_{\p_1}, R_{\p_2}, R_{\s}); ~ x_{\s}, \Sig_{\p_1}+\Sig_{\s} \Big) 
\end{align}
for a fixed given configuration $x_{\te}$ directly, we rather use the approximation 
upper bound given in \cite{patera2002method}:
\begin{align}\label{eq:algorithm_PoC_upper_bound}
    \pc\big( \cB(0, R_{\p_1}+ R_{\s}); ~ x_{\s}, \Sig_{\p_1}+\Sig_{\s} \big) + \sum_{k=1}^{N} \mathcal{P}_k
    +\pc\big( \cB(x_{\te}(1), R_{\p_2}+ R_{\s}); ~ x_{\s}, \Sig_{\p_1}+\Sig_{\s}\big),
\end{align}
where \( \mathcal{P}_k \) denotes the Gaussian integral over the swept rectangle of segment \( k \), accounting for a fixed hard-body radius.
We solve this problem using a two-phase optimization strategy:
\begin{enumerate}
    \item {Initial sampling}: multiple initial guesses are drawn at random, and configurations yielding high objective values are retained.
    \item {Gradient-based refinement}: the best candidates are refined using the constrained optimization algorithm L-BFGS-B~\cite{zhu1997algorithm_L_BFGS_B}, as implemented in the \texttt{minimize} function of the SciPy library~\cite{2020SciPy_L_BFGS_B}.

\end{enumerate}
This approach ensures both physical plausibility of the tether configurations and robustness of the risk estimate.

\end{document}